\documentclass[preprint,12pt]{elsarticle}

\makeatletter
\def\ps@pprintTitle{%
 \let\@oddhead\@empty
 \let\@evenhead\@empty
 \def\@oddfoot{\centerline{\thepage}}%
 \let\@evenfoot\@oddfoot}
\makeatother

\usepackage[latin1]{inputenc}  
\usepackage[T1]{fontenc}
\usepackage{amsmath,amssymb,dcolumn}
\usepackage{subfig}
\usepackage{graphicx,color}
\usepackage{algorithm}
\usepackage{algorithmic}
\graphicspath{{./figs/}}
\usepackage{float}
\usepackage{axecontrol,axemath,axeopt,axeabrv,axeprog,axeelsarticle}

\usepackage{ragged2e}

\begin{document}

\begin{frontmatter}
  \title{Controlling the level of sparsity in MPC} 

\author[LiU]{Daniel Axehill\corref{cor1}}
\ead{daniel@isy.liu.se}
\cortext[cor1]{Corresponding author. Tel: +46 13 284 042. Fax: +46 13
  282 622.}
\address[LiU]{Div. of Automatic Control, Dep. of Electrical Engineering,
  Linköping University, SE-581 83 Linköping, Sweden}
\begin{abstract}
  In optimization routines used for on-line Model Predictive Control
  (MPC), linear systems of equations are usually solved in each
  iteration. This is true both for Active Set (AS) methods as well as
  for Interior Point (IP) methods, and for linear MPC as well as for
  nonlinear MPC and hybrid MPC. The main computational effort is spent
  while solving these linear systems of equations, and hence, it is of
  greatest interest to solve them efficiently. Classically, the
  optimization problem has been formulated in either of two different
  ways. One of them leading to a sparse linear system of equations
  involving relatively many variables to solve in each iteration and
  the other one leading to a dense linear system of equations
  involving relatively few variables. In this work, it is shown that
  it is possible not only to consider these two distinct choices of
  formulations. Instead it is shown that it is possible to create an
  entire family of formulations with different levels of sparsity and
  number of variables, and that this extra degree of freedom can be
  exploited to get even better performance with the software and
  hardware at hand. This result also provides a better answer to an
  often discussed question in MPC; should the sparse or dense
  formulation be used. In this work, it is shown that the answer to
  this question is that often none of these classical choices is the
  best choice, and that a better choice with a different level of
  sparsity actually can be found.
\end{abstract}
\end{frontmatter}

\section{Introduction}
\label{sec:intro}
Model Predictive Control (MPC) is one of the most commonly used
control strategies in industry. Some important reasons for its success
include that it can handle multi-variable systems and constraints on
control signals and states in a structured way. In each sample, some
kind of optimization problem is solved. In the methods considered in
this paper, the optimization problem is assumed to be solved
on-line. The optimization problem can be of different types depending
on which type of system and problem formulation that is used. The most
common variants are linear MPC, non-linear MPC and hybrid MPC. In most
cases, the effort spent in the optimization problems boils down to
solving Newton-system-like equations. Hence, lots of research has been
done in the area of solving this type of system of equations
efficiently when it has the special form from MPC. It is well-known
that these equations (or at least a large part of them) can be cast in
the form of a finite horizon LQ control problem and as such it can be
solved using a Riccati recursion. Some examples of how Riccati
recursions have been used to speed up optimization routines can be
found in, for example,
\cite{jonson83:thesis,rao98:_applic_inter_point_method_model_predic_contr,hansson00:_primal_dual_inter_point_method,bartlett02:_quadr,vandenberghe02:_robus_full,akerblad04:_effic,axehill06:_mixed_integ_dual_quadr_progr_tail_mpc,axevanhan07:_relax_applic_mipc_compa_and_eff_compu,axehill08:thesis,axehill08:_dual_gradien_projec_quadr_progr,diehl09:_nonlin_model_predic_contr}.

The objective with this paper is to revisit the recurring question in
MPC whether the optimization problem should be formulated in a way
where the states are present as optimization variables or in a form
where only the control signals are the optimization variables. The
latter is often called condensing. In general this choice in turn
affects which type of linear algebra that is possible to use in the
optimization routine. If the states are kept, the system of equations
solved in each iteration in the solver potentially becomes sparse and
if the condensed formulation is used this system of equations instead
becomes dense. Since it is known that the computational complexity of
the sparse formulation grows with the prediction horizon length $N$ as
$\Ordo(N)$ if sparsity is exploited while the computational complexity
for the dense one grows as $\Ordo(N^3)$, the sparse one is often
recommended for problems with large $N$ and the dense one is often
recommended for problems with small $N$. In this paper, it will be
shown that this choice does not have to be this binary. It is shown
that it is possible to construct equivalent problems that have a level
of sparsity in between these two classical choices, but also to
increase the sparsity even further for certain types of problems. By
using the approaches proposed in this work, formulations that are even
more computationally efficient than the classical ones can be
constructed.

The key to the performance improvement is that the proposed
reformulations change the block-size of the sparse formulation. This
makes it possible to tune this choice according the performance of the
software and hardware available on-line. If a code generated solver is
used, the selection of the most appropriate formulation of the problem
can be done off-line in the code generation phase, taking the
performance of the on-line hardware and software platform into
account. Choosing a good, or even ``optimal'', block-size is not a new
idea in numerical linear algebra. See, \eg,
\cite{whaley05:_minim_blas}. However, it has, to the best of our
knowledge, not previously been discussed for MPC and in particular not
in the context of tailored linear algebra for MPC. Furthermore, the
presented result hints that control engineers working with MPC should
not ask themselves whether to formulate the problem in a sparse or
dense way, but instead what is the correct \emph{level} of sparsity
for the problem at hand to obtain maximum performance which do not
necessarily coincide with one of the two extreme choices that have
been used classically.

In this article, $\symmats^n$ denotes symmetric matrices with $n$
columns. Furthermore, $\posdefmats^n$ ($\possemidefmats^n$) denotes
symmetric positive (semi) definite matrices with $n$ columns. The set
$\posintnums$ denote the set of positive non-zero integers. A Sans
Serif font is used to indicate that a matrix or a vector is, in some
way, composed of stacked matrices or vectors from different time
instants. The stacked matrices or vectors have a similar symbol as the
composed matrix, but in an ordinary font. For example, $\timestack{Q}
= \diag(Q,\hdots,Q)$.

\section{Problem formulation}
The problem considered in this work is
\begin{equation}
  \label{eq:lq_opt_prob_def}
  \minimize{\sum_{t=0}^{N-1}\frac{1}{2}
    \begin{bmatrix}
      x_t^T & u_t^T
    \end{bmatrix}
    \begin{bmatrix}
      Q_t & W_t\\
      W_t^T & R_t
    \end{bmatrix}
    \begin{bmatrix}
      x_t \\
      u_t
    \end{bmatrix}
    + \frac{1}{2}x_N^TQ_Nx_N}
  {x_t,u_t}
  {x_0 &= x^0\\
  x_{t+1} &= A_tx_t + B_tu_t,\; \timespan{0}{N-1}\\
  0 &\geq H_{x,t}x_t+H_{u,t}u_t + h_t,\; \timespan{0}{N-1}\\
0 &\geq H_{x,N}x_N + h_N}
\end{equation}
where the states $x_t \in \realnums^{n}$, the initial condition $x^0
\in \realnums^{n}$, the control inputs $u_t \in \realnums^{m}$, the
system matrices $A_t \in \realnums^{n\times n}, B_t \in \realnums^{n
  \times m}$, the penalty matrices for the states $Q_t \in
\symmats^n$, penalty matrices for the control inputs $R_t \in
\symmats^m$, the cross penalty matrices $W_t \in \realnums^{n \times
  m}$, the state constraint coefficient matrix $H_{x,t} \in
\realnums^{c \times n}$, the control signal coefficient matrix
$H_{u,t} \in \realnums^{c \times m}$, the constraint constant $h_t \in
\realnums^c$, and the prediction horizon $N \in
\nonnegintnums$. Moreover, the matrices $Q_t$, $R_t$ and $W_t$ are
assumed to be chosen such that the following two assumptions are
satisfied
\begin{ass}
  \label{ass:R_t_pos_def}
  $R_t \in \posdefmats^m$
\end{ass}
\begin{ass}
  \label{ass:penalty_mat_pos_def}
  $\begin{bmatrix} Q_t & W_t\\
    W_t^T & R_t \end{bmatrix}
  \in \possemidefmats^{n+m}$
\end{ass}
Both constrained linear MPC problems and nonlinear MPC problems often
boil down to solving problems similar to the one
in~\eqref{eq:lq_opt_prob_def} but without any inequality constraints
during the Interior Point (IP) process or Active Set (AS) process,
\cite{jonson83:thesis,rao98:_applic_inter_point_method_model_predic_contr,hansson00:_primal_dual_inter_point_method,bartlett02:_quadr,vandenberghe02:_robus_full,akerblad04:_effic,axehill06:_mixed_integ_dual_quadr_progr_tail_mpc,axevanhan07:_relax_applic_mipc_compa_and_eff_compu,axehill08:thesis,axehill08:_dual_gradien_projec_quadr_progr,diehl09:_nonlin_model_predic_contr}. Hence,
the ability of solving unconstrained versions of the problem
in~\eqref{eq:lq_opt_prob_def} efficiently is of great interest for the
overall computational performance in the entire range of problems from
simple unconstrained linear MPC problems, to nonlinear constrained and
hybrid MPC problems.

As shown in~\cite{astrom84:_comput_contr_system}, the problem
in~\eqref{eq:lq_opt_prob_def} can after a simple variable
transformation be recast in an equivalent form with
$W_t=0$. Therefore, the analysis in this work is restricted to the
case when $W_t=0$ without any loss of generality.

\begin{ass}
  \label{ass:cross_penalty_zero}
  $W_t=0$
\end{ass}

\begin{rmk}
  The results shown in this paper can easily be extended to common
  variants of MPC, \eg, problems where the control signal horizon
  differs from the prediction horizon, problems with affine system
  descriptions, as well as to reference tracking problems.
\end{rmk}

\section{Classical optimization problem formulations of MPC}
\label{sec:class_qp_forms}
Traditionally, two optimization problem formulations of the MPC
problem have been dominating in the MPC community. In the first
formulation, the optimization problem in~\eqref{eq:lq_opt_prob_def}
with $W_t=0$ has been written more compactly as a Quadratic
Programming (QP) problem in the form
\begin{equation}
  \label{eq:qp_sparse_form}
  \minimize{\frac{1}{2}
    \begin{bmatrix}
      \timestack{x}^T & \timestack{u}^T
    \end{bmatrix}
    \begin{bmatrix}
      \timestack{Q} & 0 \\
      0 & \timestack{R}
    \end{bmatrix}
    \begin{bmatrix}
      \timestack{x} \\
      \timestack{u}
    \end{bmatrix}}
  {\timestack{x},\timestack{u}}
  {\begin{bmatrix}
     \timestack{A} & \timestack{B}
   \end{bmatrix}
   \begin{bmatrix}
      \timestack{x} \\
      \timestack{u}
    \end{bmatrix}
    &= \timestack{b}\\
  \begin{bmatrix}
     \timestack{H}_{\timestack{x}} & \timestack{H}_{\timestack{u}}
   \end{bmatrix}
   \begin{bmatrix}
      \timestack{x} \\
      \timestack{u}
    \end{bmatrix}
    + \timestack{h}
    &\leq 0}
\end{equation}
where $\timestack{x}$, $\timestack{u}$, $\timestack{Q}$,
$\timestack{R}$, $\timestack{q}$, $\timestack{r}$, $\timestack{b}$,
$\timestack{A}$, $\timestack{B}$, $\timestack{H}_{\timestack{x}}$,
$\timestack{H}_{\timestack{u}}$, and $\timestack{h}$ are defined in
the Appendix. Note that, $\timestack{Q}$, $\timestack{R}$,
$\timestack{A}$, $\timestack{B}$, $\timestack{H}_{\timestack{x}}$, and
$\timestack{H}_{\timestack{u}}$ are sparse matrices.

In the second formulation, the dynamics equations
in~\eqref{eq:lq_opt_prob_def} have been used to express
$\timestack{x}$ as
\begin{equation}
  \label{eq:MPC_X_of_U_and_x0}
  \timestack{x} = S_{x,N} x_0+S_{u,N} \timestack{u}
\end{equation}
where $S_{x,N}$ and $S_{u,N}$ are defined in the Appendix. This
expression can be used to eliminate the equality constraints
containing the dynamics in the problem
in~\eqref{eq:qp_sparse_form}. As a result, an equivalent problem can
be derived in the form
\begin{equation}
  \label{eq:qp_dense_form}
  \minimize{&\frac{1}{2}\timestack{u}^T\parens{\timestack{B}^T\timestack{A}^{-T}\timestack{Q}\timestack{A}^{-1}\timestack{B}
      + \timestack{R}}\timestack{u}-\parens{\timestack{B}^T\timestack{A}^{-T}\timestack{Q}\timestack{A}^{-1}\timestack{b}}^T\timestack{u}}
  {\timestack{u}}
  {\parens{\timestack{H}_{\timestack{x}}S_{u,N} +
      \timestack{H}_{\timestack{u}}}\timestack{u} + \timestack{h} +
    \timestack{H}_{\timestack{x}}S_{x,N}x_0 \leq 0}
\end{equation}
where
$\timestack{B}^T\timestack{A}^{-T}\timestack{Q}\timestack{A}^{-1}\timestack{B}
+ \timestack{R}$ is a dense matrix.

\section{Quasi-sparse optimization formulations}
As discussed in the introduction, many papers published on the subject
illustrate that the sparse formulation in~\eqref{eq:qp_sparse_form} is
preferable over the dense one in~\eqref{eq:qp_dense_form} from a
computational point of view. However, there also exists applications
where the non-structure exploiting dense formulation turns out to be
the fastest one. It can be realized both from the expressions for the
analytical complexities as well as from numerical experiments that
there are certain breakpoints in the problem sizes where one
formulation is better than the other one. Traditionally, one
rule-of-thumb is that the sparse formulation is faster for large
values of $N$, and the dense one is faster for problems with small
values of $N$. In this section it is discussed whether it is possible
to do even better by combining ideas from these two traditional
approaches. To reduce the complexity of the notation in the
presentation, the ideas are illustrated on an MPC formulation where
the system, penalties and constraints are independent of time.

\subsection{Increasing the block-size}
\label{subsec:inc_blk_size}
Partition the prediction horizon in $\tilde{N}+1$ subintervals with
corresponding lengths $M_k \in \posintnums$, $k \in \{0,\tilde{N}\}$
with $\sum_{k=0}^{\tilde{N}} M_k = N+1$. To get a reformulated problem
which is in the form of the one in~\eqref{eq:lq_opt_prob_def} with a
non-zero end penalty, the choice $M_{\tilde{N}}=1$ was made in this
section. Define $\timestack{x}_k = [x_{\tau_k+0}^T \, x_{\tau_k+1}^T
\hdots x_{\tau_k+M_k}^T]^T$ and $\timestack{u}_k = [u_{\tau_k+0}^T \,
u_{\tau_k+1}^T \hdots u_{\tau_k+M_k-1}^T]^T$ where $\tau_0=0$ and
$\tau_k = \sum_{i=0}^{k-1} M_i, \; k > 0$. Analogously to the equation
in~\eqref{eq:MPC_X_of_U_and_x0}, given the state at a time $\tau_k$
and all control signals from subinterval $k$, all states in
subinterval $k$ can be expressed as
\begin{equation}
  \label{eq:MPC_X_of_U_and_x_tau}
  \timestack{x}_k = S_{x,M_k} x_{\tau_k}+S_{u,M_k} \timestack{u}_k
\end{equation}
and in particular we have that
\begin{equation}
  \label{eq:MPC_x_tau_p_N_of_U_and_x_tau}
  x_{\tau_{k+1}} = x_{\tau_k+M_k} = A^{M_k} x_{\tau_k}+
  \begin{bmatrix}
    A^{M_k-1}B & A^{M_k-2}B & \hdots & B
  \end{bmatrix}\timestack{u}_k
  \triangleq \timestack{A}_kx_{\tau_k}+\timestack{B}_k\timestack{u}_k
\end{equation}
Note that the equation in~\eqref{eq:MPC_x_tau_p_N_of_U_and_x_tau} is
in state-space form. The new state dynamics describes the dynamics
from the first sample in one condensed block to the first sample in
the following block. Using these results, it is possible to write the
sum of the stage costs for an entire subinterval as
\begin{equation}
  \label{eq:substage_cost}
  \sum_{t=\tau_k}^{\tau_k+M_k-1}
  x_{t+1}^TQx_{t+1} + u_t^TRu_t=
  \begin{bmatrix}
      x_{\tau_k}^T & \timestack{u}_k^T
    \end{bmatrix}
    \begin{bmatrix}
      \timestack{Q}_k & \timestack{W}_k\\
      \timestack{W}_k^T & \timestack{R}_k
    \end{bmatrix}
    \begin{bmatrix}
      x_{\tau_k} \\
      \timestack{u}_k
    \end{bmatrix}
\end{equation}
with
\begin{equation}
  \label{eq:substage_matrices_def}
  \begin{split}
    \timestack{Q}_k &=
    S_{x,M_k}^T\cdot\diag\parens{Q,\hdots,Q,
      0}\cdot S_{x,M_k},\\
    \timestack{W}_k &=
    S_{x,M_k}^T\cdot\diag\parens{Q,\hdots,Q,
      0}\cdot
    S_{u,M_k},\\
    \timestack{R}_k &=
    S_{u,M_k}^T\cdot\diag\parens{Q,\hdots,Q,
      0}\cdot
    S_{u,M_k}+\diag\parens{R,\hdots,R}
  \end{split}
\end{equation}
for $k \leq \tilde{N}-1$ and $\timestack{Q}_{\tilde{N}} = Q$.  As a
result, the original problem can be re-cast in the form
\begin{equation}
  \label{eq:lq_opt_prob_def_mergestage}
  \minimize{\sum_{k=0}^{\tilde{N}-1}\frac{1}{2}
    \begin{bmatrix}
      x_{\tau_k}^T & \timestack{u}_k^T
    \end{bmatrix}
    \begin{bmatrix}
      \timestack{Q}_k & \timestack{W}_k\\
      \timestack{W}_k^T & \timestack{R}_k
    \end{bmatrix}
    \begin{bmatrix}
      x_{\tau_k} \\
      \timestack{u}_k
    \end{bmatrix}
  + \frac{1}{2}x_{\tau_{\tilde{N}}}^T \timestack{Q}_{\tilde{N}}x_{\tau_{\tilde{N}}}}
  {x_{\tau_k},\timestack{u}_k}
  {x_{\tau_0} &= x^0\\
  x_{\tau_{k+1}} &= \timestack{A}_kx_{\tau_{k}} +
  \timestack{B}_k\timestack{u}_k,\; \timespan[k]{0}{\tilde{N}-1}\\
  0 &\geq \timestack{H}_{x,k}x_{\tau_k} +
  \timestack{H}_{u,k}\timestack{u}_k + \timestack{h}_{k},\;
  \timespan[k]{0}{\tilde{N}-1}\\
  0 &\geq \timestack{H}_{x,\tilde{N}}x_{\tau_{\tilde{N}}} + \timestack{h}_{\tilde{N}}}
\end{equation}
with
\begin{equation}
  \label{eq:constr_mat_merge_def}
  \begin{split}
  \timestack{H}_{x,k} &=
  \begin{bmatrix}
    H_{x} & 0 & \hdots & 0 & 0\\
    0 & \ddots & \hdots & 0 & 0\\
    0 & \hdots & 0 & H_{x} & 0
  \end{bmatrix}
    S_{x,M_k} \in \realnums^{M_kc \times (M_k+1)n},\\
  \timestack{H}_{u,k} &=
  \begin{bmatrix}
    H_{x} & 0 & \hdots & 0\\
    0 & \ddots & \hdots & 0\\
    0 & \hdots & H_{x} & 0
  \end{bmatrix}S_{u,M_k} + 
  \begin{bmatrix}
    H_{u} & 0 & \hdots & 0\\
    0 & \ddots & \hdots & 0\\
    0 & \hdots & 0 & H_{u}
  \end{bmatrix} \in \realnums^{M_kc \times M_km},\\
  \timestack{h}_{k} &= \begin{bmatrix}h \\ \vdots \\
    h\end{bmatrix}\in \realnums^{M_kc},\; \timespan[k]{0}{\tilde{N}-1}\\
  \timestack{H}_{x,\tilde{N}} &= H_{x}, \; \timestack{h}_{\tilde{N}} = h
\end{split}
\end{equation}
The new formulation can be interpreted as another MPC problem in the
form in~\eqref{eq:lq_opt_prob_def} with virtual prediction horizon
$\tilde{N}$, virtual state dimension $\tilde{n} = n$, and the virtual
control signal dimension for interval $k$ is $\tilde{m} = M_k \cdot
m$. There are different variants of this formulation that give similar
results. For example one can take $\timestack{Q}_{\tilde{N}}=0$ and
instead include the last state in the second last sub-interval. The
number of inequality constraints in each virtual time step along the
prediction horizon is $\tilde{c}_k = M_k \cdot c$ and the total number
of inequality constraints is unaffected compared to the original
problem. In words, this formulation partially condenses the original
sparse problem into a new one where several original samples have been
condensed into one new and are basically handled using the dense
formulation. If $M_0=N$, roughly the traditional dense formulation
in~\eqref{eq:qp_dense_form} is obtained and if $M_k=1,\;\forall k$,
the traditional sparse one in~\eqref{eq:qp_sparse_form} is
obtained. Note that, the problems in~\eqref{eq:lq_opt_prob_def} and
in~\eqref{eq:lq_opt_prob_def_mergestage} are structurally identical
and an algorithm (and an implementation of it) that can be applied to
the problem in~\eqref{eq:lq_opt_prob_def} can also be applied to the
one in~\eqref{eq:lq_opt_prob_def_mergestage}. Once the problem
in~\eqref{eq:lq_opt_prob_def_mergestage} has been solved,
$\timestack{u}$ in Section~\ref{sec:class_qp_forms} is directly
obtained and the entire vector $\timestack{x}$ can, if desired, easily
be computed.

\subsection{Decreasing the block-size}
\label{subsec:dec_blk_size}
Consider a problem with diagonal $R$ and inequality constraints in
the form
\begin{equation}
  \label{eq:state_upd_ineq_constr}
  0 \geq
  \begin{bmatrix}
    H_{x} & 0\\
    0 & H_{u}
  \end{bmatrix}
  \begin{bmatrix}
    x(t)\\
    u(t)
  \end{bmatrix}
   +
   \begin{bmatrix}
     h_{x}\\
     h_{u}
   \end{bmatrix}
\end{equation}
where $H_{u}$ represents simple control signal constraints. Simple
here means that the control signal constraints are, \eg, upper and
lower bound constraints. For simplicity, we assume $H_{u}$ diagonal
with diagonal elements $H_{u,k}, k = 1,\hdots,m$ in what
follows. Note that, however, $H_{x}$ can be full as long as these
state constraints do not involve control signals.

Partition the control signal $u_t$ in $M_t$ parts such that $u_t =
\begin{bmatrix}u_{t,1} & u_{t,2} & \hdots & u_{t,M_t}\end{bmatrix}$
with $u_{t,k} \in \realnums^{m,\tilde{m}_{t,k}}$ and $B =
\begin{bmatrix}b_1 & b_2 & \hdots & b_{M_t}\end{bmatrix}$, where
$\tilde{m}_{t,k} \in \posintnums$, $t \in \{0,N-1\}$, $k \in
\{1,M_t\}$ denotes the number of virtual control signals for which it
holds that $\sum_{t=0}^{N-1}\sum_{k=1}^{M_t} \tilde{m}_{t,k} =
Nm$. Then the state update equation can be written as
\begin{equation}
  \label{eq:state_upd_eq_sum}
  x_{t+1} = Ax_t + Bu_t = Ax_t + \sum_{k=1}^{M_t}b_ku_{t,k}
\end{equation}
which can be equivalently reformulated, \eg, as
\begin{equation}
  \label{eq:state_upd_eq_substage}
  \begin{split}
    \tilde{x}_{t,2} &= A\tilde{x}_{t,1} + b_{1}u_{t,1}\\
    \tilde{x}_{t,3} &= I\tilde{x}_{t,2} + b_{2}u_{t,2}\\
    &\vdots\\
    \tilde{x}_{t+1,1} &= I\tilde{x}_{t,M_t} + b_{M_t}u_{t,M_t}
  \end{split}
\end{equation}
with some new states $\tilde{x}_{t,k}$ for which it holds that
$\tilde{x}_{t,1}=x_{t}, \; \forall t$. Since $R$ and $H_{u}$ are
diagonal, the objective function and inequality constraints can be
decomposed in $u(t)$ in an analogous way and the diagonal matrices
corresponding to part $k$ are denoted $R_k$ and $H_{u,k}$,
respectively. That in combination with the expression
in~\eqref{eq:state_upd_eq_substage}, the problem
in~\eqref{eq:lq_opt_prob_def} can be reformulated in the form
\begin{equation}
  \label{eq:lq_opt_prob_def_substage}
  \minimize{\sum_{k=0}^{\tilde{N}-1}\frac{1}{2}
    \begin{bmatrix}
      \tilde{x}_{k}^T & \timestack{u}_k^T
    \end{bmatrix}
    \begin{bmatrix}
      \timestack{Q}_k & \timestack{W}_k\\
      \timestack{W}_k^T & \timestack{R}_k
    \end{bmatrix}
    \begin{bmatrix}
      \tilde{x}_{k} \\
      \timestack{u}_k
    \end{bmatrix}
  + \frac{1}{2}\tilde{x}_{\tilde{N}}^T \timestack{Q}_{\tilde{N}}\tilde{x}_{\tilde{N}}}
  {\tilde{x}_{k},\timestack{u}_k}
  {\tilde{x}_{0} &= x^0\\
  \tilde{x}_{k+1} &= \timestack{A}_k\tilde{x}_{k} +
  \timestack{B}_k\timestack{u}_k,\; \timespan[k]{0}{\tilde{N}-1}\\
  0 &\geq \timestack{H}_{x,k}\tilde{x}_{k} +
  \timestack{H}_{u,k}\timestack{u}_k + \timestack{h}_{k},\;
  \timespan[k]{0}{\tilde{N}-1}\\
  0 &\geq \timestack{H}_{x,\tilde{N}}\tilde{x}_{\tilde{N}} + \timestack{h}_{\tilde{N}}}
\end{equation}
 with
\begin{equation}
  \label{eq:substage_def}
  \begin{split}
    \timestack{A}_k &=
    \begin{cases}
      &A, \; k=0,M,2M,\hdots,(N-1)M\\
      &I, \; \text{otherwise}
    \end{cases}\\
    \timestack{B}_k &= b_{\text{mod}(k,M)+1}\\
    \timestack{Q}_k &=
    \begin{cases}
      &Q, \; k=0,M,2M,\hdots, NM\\
      &0, \; \text{otherwise}
    \end{cases}\\
    \timestack{R}_k &= R_{\text{mod}(k,M)+1}\\
    \timestack{H}_{x,k} &=
    \begin{cases}
      &
      \begin{bmatrix}
        H_{x}\\
        0
      \end{bmatrix}, \; k=0,M,2M,\hdots,(N-1)M\\
      &H_{x}, \; k=\tilde{N}\\
      &0, \; \text{otherwise}
    \end{cases}\\
    \timestack{H}_{u,k} &= 
    \begin{cases}
      &
      \begin{bmatrix}
        0\\
        H_{u,1}
      \end{bmatrix}, \; k=0,M,2M,\hdots,(N-1)M\\
      &H_{u,\text{mod}(k,M)+1}, \text{otherwise}
    \end{cases}\\
    \timestack{h}_k &=
    \begin{cases}
      &
      \begin{bmatrix}
        h_x\\
        h_{u,1}
      \end{bmatrix}, \; k=0,M,2M,\hdots,(N-1)M\\
      &h_x, \; k=\tilde{N}\\
      &h_{u,\text{mod}(k,M)+1}, \; \text{otherwise}
    \end{cases}
  \end{split}
\end{equation}
for the simplified case when the partitioning is done uniformly with a
constant $M_t = M \in \braces{j \in\posintnums: m/j\in\posintnums},\;
\forall t$ over the prediction horizon. The reformulated problem can
be interpreted as an optimization problem in the MPC form
in~\eqref{eq:lq_opt_prob_def} with virtual prediction horizon
$\tilde{N} = NM$, virtual state dimension $\tilde{n} = n$, and virtual
control signal dimension $\tilde{m} = m/M$. The number of inequality
constraints is varying along the virtual prediction horizon, however,
the total number of constraints is unchanged compared to the original
problem. Once the problem in~\eqref{eq:lq_opt_prob_def_substage} has
been solved, $\timestack{u}$ and $\timestack{x}$ in
Section~\ref{sec:class_qp_forms} are directly obtained.

\begin{rmk}
  \label{rmk:combination}
  Further degrees of freedom can be obtained by combining the
  approaches in Sections \ref{subsec:inc_blk_size} and
  \ref{subsec:dec_blk_size} by first decreasing the block-size and
  then increasing the block-size of the modified problem.
\end{rmk}

\begin{rmk}
  \label{rmk:filtering}
  Note that, the ideas can also be applied to problems with similar
  structures. For example, to the moving horizon state-estimation
  problem and to the Lagrange dual of the optimal control problem.
\end{rmk}

\section{Impact on a commonly used sparse linear algebra for MPC}
\label{sec:theo_impact_flops}
The inequality constrained optimal control problems
in~\eqref{eq:lq_opt_prob_def_mergestage} and
\eqref{eq:lq_opt_prob_def_substage} are usually either solved using
either an interior point (IP) method or an active-set (AS) method. It
is well-known that the main computational complexity in these
algorithms can be formulated as solving a sequence of unconstrained
variants of the original control problem. These problems in turn are
solved by solving a linear system of equations corresponding to the
KKT conditions of these unconstrained problems. This can be done in
several ways. However, two commonly used approaches are to either use
a Cholesky factorization applied to a problem where the states have
been eliminated (commonly known as condensing) or a Riccati
factorization applied to the problem where the states are kept as
variables. If the problem is reformulated as described in
Sections~\ref{subsec:inc_blk_size} and \ref{subsec:dec_blk_size}, it
can be realized that the important block-sizes that appear in a sparse
factorization will be changed. Even thought the Riccati factorization
is used as an example in this work, it is expected that the proposed
approaches in this work will have similar impact on other variants of
sparse linear algebra used in MPC. However, none of the proposed
reformulations will affect the dense formulation and the following
(off-the-shelf) dense linear algebra. Hence, the focus in this section
will be on how Riccati based sparse linear algebra is affected by the
approaches introduced in this article. The required number of flops
for the Cholesky factorization approach is known to be roughly
$(Nm)^3/3$ and for the Riccati factorization it is roughly
$N(1/3m^3+4n^3+4m^2n+6n^2m)$,
\cite{nielsen13:_low_modif_riccat_factor_applic_unpub}. A full
specification of the Riccati approach used in this work can be found
in~\cite{nielsen13:_low_modif_riccat_factor_applic_unpub,axehill08:thesis}. For
simplicity, in this section it is assumed that the blocking is uniform
non-time-dependent and that it is compatible with $N$ and $m$. In
practice further improvements can potentially be achieved by using a
non-uniform blocking factor. However, this also increases the
complexity of determining a good one and it is therefore not clear
whether it is worth that extra effort. Furthermore, it would also be
possible to exploit and take into account in the complexity
calculations the special structures in the reformulated problems. For
example, the approach in Section~\ref{subsec:dec_blk_size} generates
many sparse problem data matrices.

\subsection{Increasing the block-size}
\label{subsec:inc_blk_size_impact}
\begin{figure}[tbh]
  \centering
  \includegraphics[width=0.45\textwidth]{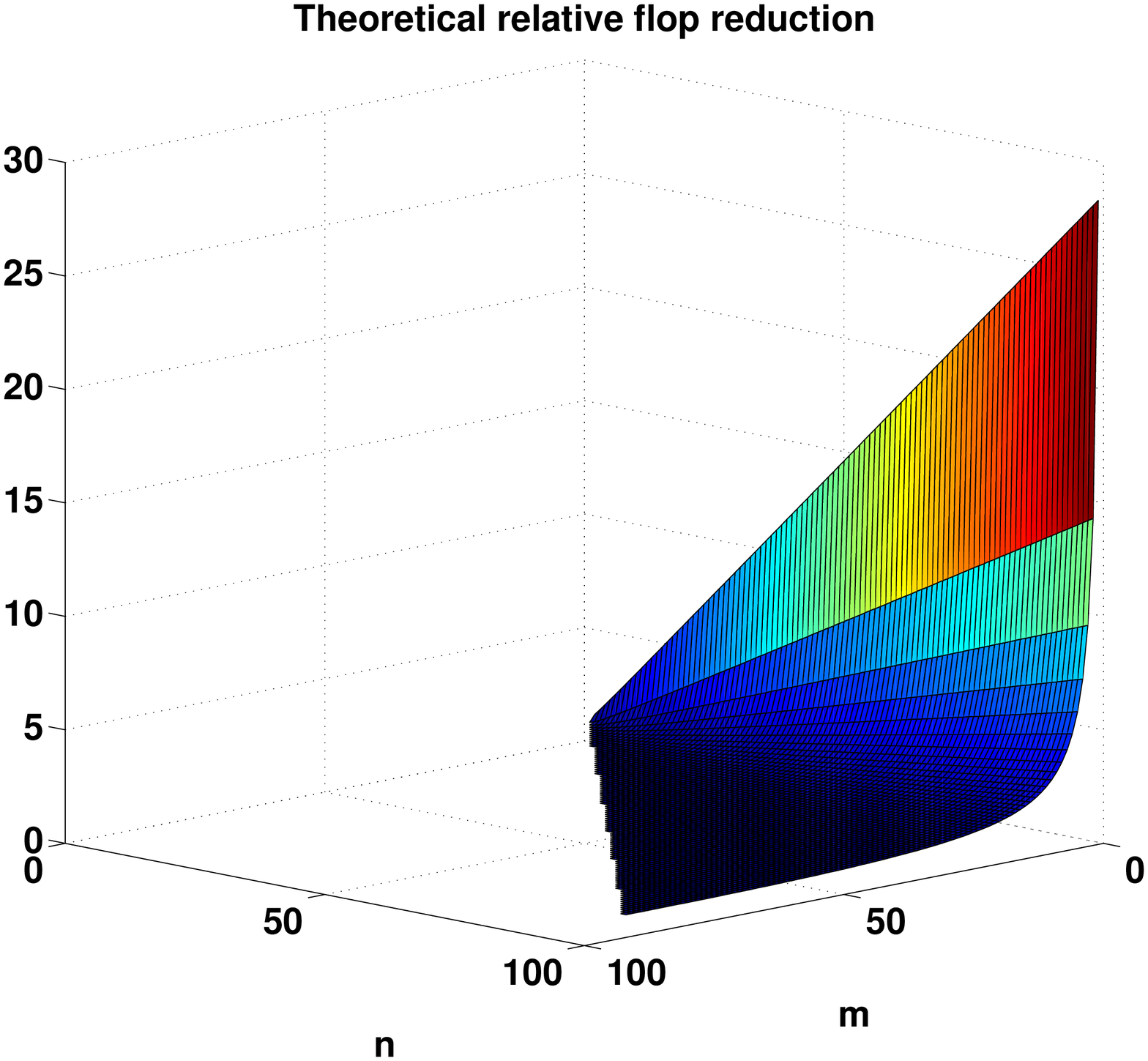}\hspace{0.1\textwidth}\includegraphics[width=0.45\textwidth]{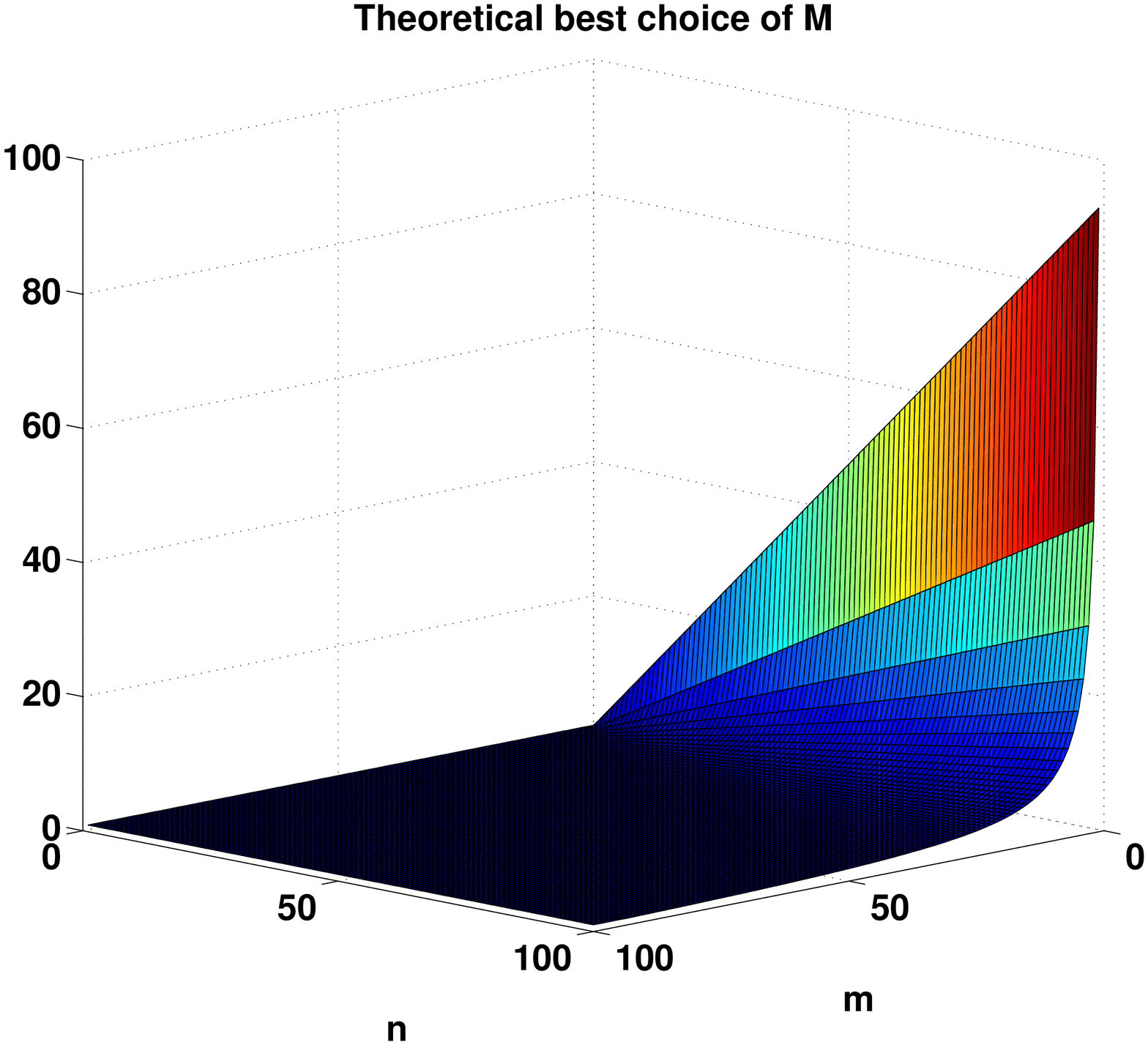}
  \caption{Theoretical maximum relative flop reduction (left) and for
    which respective choice of $M$ this is obtained (right) for a
    range of relevant values of $n$ and $m$. In this range a flop
    reduction of almost up to $30$ times can be achieved. The
    reduction in the left plot is only plotted for combinations of $n$
    and $m$ that will give an improvement.}
  \label{fig:appr_1_rel_flop_gain}
\end{figure}

The family of formulations achievable with this approach is
parameterized by $M$, which will increase the block-size in the
resulting sparse numerical linear algebra and make the formulation
less sparse as $M$ grows. Since it holds that $\tilde{N}=N/M$,
$\tilde{n} = n$, and $\tilde{m}=Mm$, the required number of flops for
performing a Riccati factorization on the reformulated problem is
$f_1(M)=N\parens{M^2m^3/3 + 4n^3/M + 4Mm^2n + 6n^2m}$. The theoretical
maximum possible gain in flops $f_1(1)/f_1(M^*)$ which is independent
of $N$ is illustrated in Figure~\ref{fig:appr_1_rel_flop_gain} and can
be found to be almost as much as a factor of $30$ in the considered
problems. $M^*$ denotes the best choice of $M$ for given $n$ and
$m$. Furthermore, note that, the relative improvement $f_1(1)/f_1(M)$
is independent of $N$. However, given an $N$ one can only select an $M
\leq N$.
From the plot it follows that this approach is useful if $m$ is small
compared to $n$. In that case, it is beneficial to reformulate the
problem as an equivalent control problem with a shorter virtual
prediction horizon with more virtual control signals in each virtual
sample.

\subsection{Decreasing the block-size}
\label{subsec:dec_blk_size_impact}
\begin{figure}[tbh]
  \centering
  \includegraphics[width=0.45\textwidth]{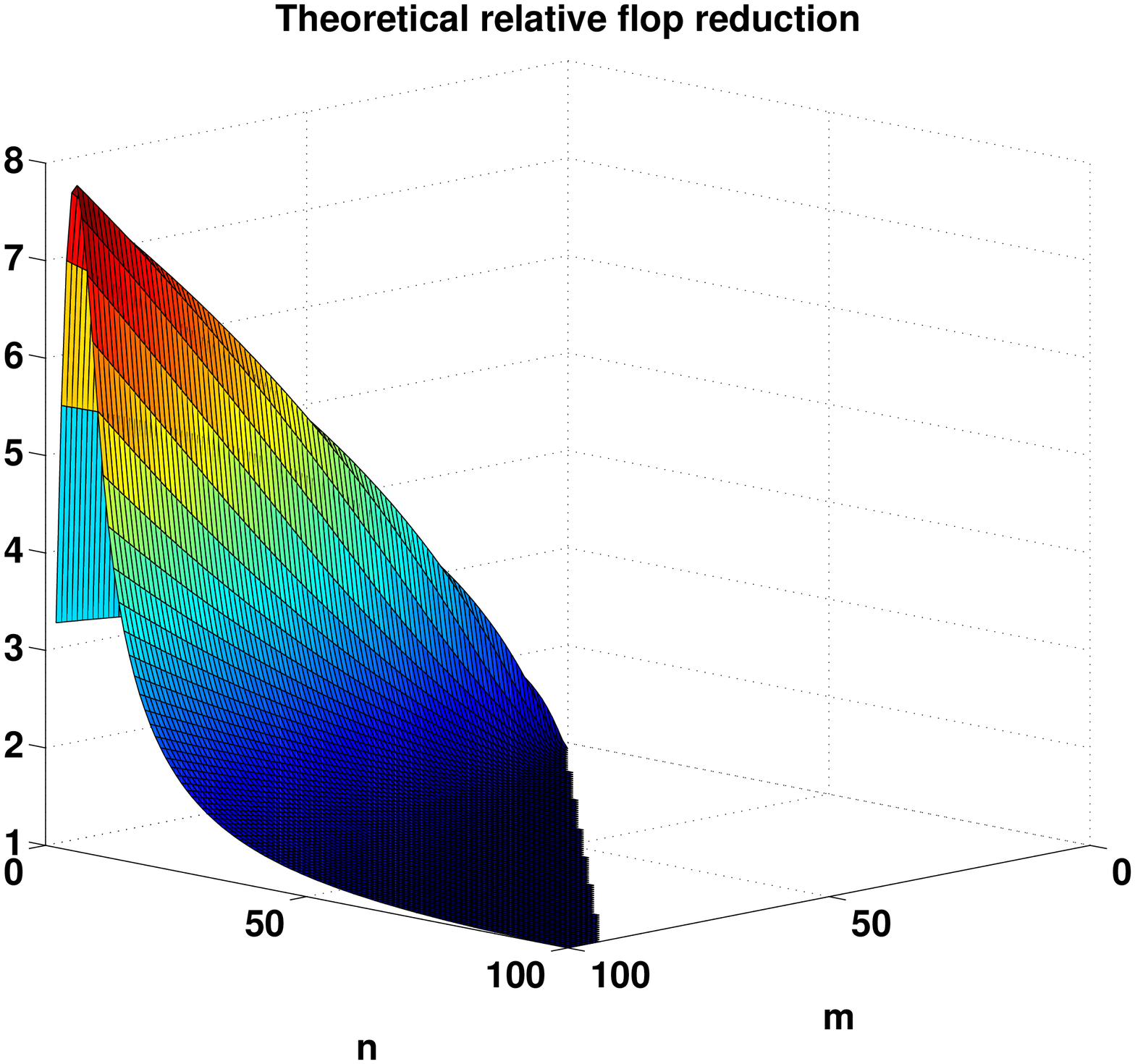}\hspace{0.1\textwidth}\includegraphics[width=0.45\textwidth]{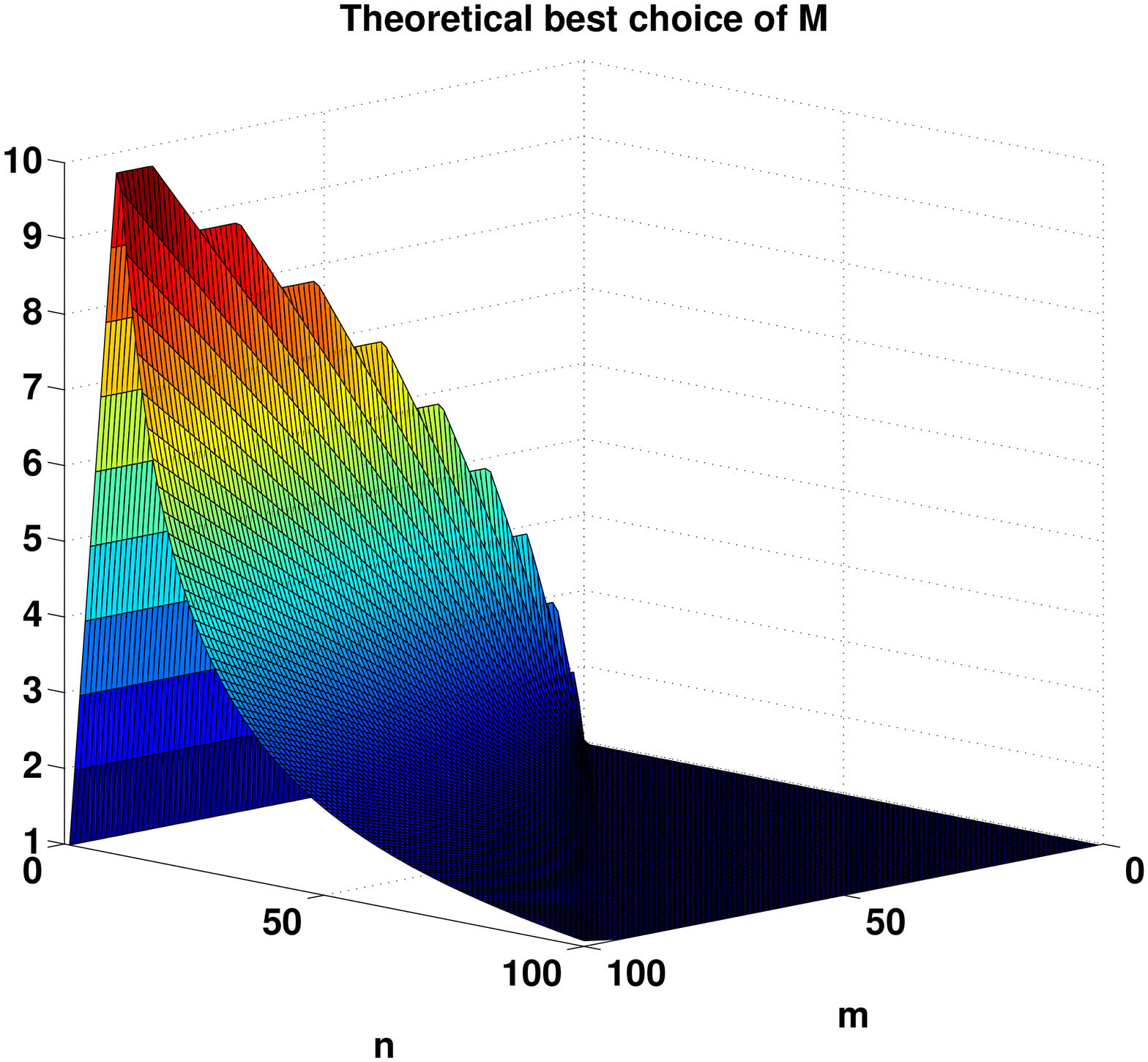}
  \caption{Theoretical maximum relative flop reduction (left) and for
    which respective choice of $M$ this is obtained (right) for a
    range of relevant values of $n$ and $m$. In this range a flop
    reduction of almost up to $8$ times can be achieved. The reduction
    in the left plot is only plotted for combinations of $n$ and $m$
    that will give an improvement.}
  \label{fig:appr_2_rel_flop_gain}
\end{figure}

This approach is related to what in the Kalman filtering literature is
known as sequential processing, where in certain cases the measurement
update can be performed sequentially for each one of the measurements
in each sample,
\cite{kailath00:_linear,gustafsson00:_adapt_filter_chang_detec}. The
family of formulations achievable with this approach is again
parameterized by $M$, which will here decrease the block-size in the
resulting sparse numerical linear algebra and make the formulation
more sparse as $M$ grows. Since $\tilde{N}=NM$, $\tilde{n} = n$, and
$\tilde{m}=m/M$ and the required number of flops for performing a
Riccati factorization on the reformulated problem is $f_2(M) =
N\parens{m^3/(3M^2) + 4n^3M + 4m^2n/M + 6n^2m}$. The theoretical
maximum possible gain in flops $f_2(1)/f_2(M^*)$ is illustrated in
Figure~\ref{fig:appr_2_rel_flop_gain} and can be found to be almost as
much as a factor of $8$ in the considered problems.
From the plot it follows that this approach is useful if $n$ is small
compared to $m$. In that case, it is beneficial to reformulate the
problem with a larger prediction horizon with less control signals in
each sample.

Note that, by combining the results from
Figures~\ref{fig:appr_1_rel_flop_gain} and
\ref{fig:appr_2_rel_flop_gain} it can be seen that where the
performance gain of one approach is lost, the other one starts to
offer a gain. Still, the improvement possible by the presented
approaches is only moderate for problems for which $n \approxeq m$,
which indicates that the sparsity obtained from the original
formulation is a good choice from a performance point of view.

\subsection{Combining the approaches}
\label{subsec:dec_blk_sizecombined_impact}
The presented approaches rely on that the original problem's
prediction horizon either is virtually reduced with the cost of having
more control signals in each virtual sample, or that the number of
control signals is reduced in each virtual sample with the cost of
getting a longer virtual prediction horizon. A natural extension is
to have a virtual sampling rate over the virtual horizon that is not
in sync with the original sampling rate. This can be achieved by
combining the two proposed procedures as pointed out in
Remark~\ref{rmk:combination} by first increasing the prediction
horizon length according to Section~\ref{subsec:dec_blk_size} and then
reducing it again according to Section~\ref{subsec:inc_blk_size} with
another choice of $M_k$ (if the choice of $M_k$ is the same in both
steps the end result would be the original problem formulation).

\begin{rmk}
  \label{rmk:parallel}
  The focus in the numerical experiments in this section is serial
  linear algebra. However, more generally, the proposed reformulations
  are also interesting for parallel approaches where they for example
  potentially can be used to optimize the workload distribution.
\end{rmk}

\section{Numerical experiments}
In these numerical examples the performance of the proposed strategies
is evaluated. The purpose is to investigate how the theoretical
improvements in terms of flops in Section~\ref{sec:theo_impact_flops}
are transfered to improvements in terms of computational time. The
experiments are performed on an Intel \mbox{i5-2520M} with $8\,$GiB
RAM running Windows~7 64-bit. All algorithms, including the Cholesky
factorizations, are implemented in m-code in an attempt to make the
results relevant for cases where everything is written in the same
programming language. If the Cholesky factorization would have been
carried out using, \eg, LAPACK calls then the proposed tools could
have been used to find other values of $M$ to maximize the overall
performance. The reformulations presented in this paper are performed
in a preprocessing stage, and the reformulated problems are sent to an
implementation of a sparse Riccati KKT system solver. As a comparison,
the same problem is also solved, first, using a standard dense
formulation using a Cholesky factorization and, second, using the
sparse Riccati KKT solver applied to the initial formulation of the
problem (block-sizes given by the problem). This experiment basically
illustrates the computational complexity of the step direction
computation in an optimization routine. To simulate the computations
for the search step direction computation in an AS or IP method, the
test problems considered are in the form in~\eqref{eq:lq_opt_prob_def}
without any inequality constraints. Note that, the actual choice of
examples used in the experiments are irrelevant since the performance
of the used linear algebra is only affected by the sizes of problem
matrices and not the numbers contained in these. Since the computed
direction from the reformulated problem is the same as from the
original formulation, apart from minor differences due to differences
in the numerics, the number of steps performed by a solver can be
expected to be the same independently of which formulation that is
used. Since usually the main computational effort in the targeted
optimization routines originates from the step direction computation,
the overall computational time can be expected to roughly scale as the
computational time for the linear algebra which is what is
investigated in this section. Furthermore, the KKT solvers used do not
at all utilize the structure in the sub-blocks of the reformulated
problems, which means that the results shown here can easily be
obtained by users without actually making any other changes of their
MPC codes rather than the problem formulation as described in earlier
sections. It can at least in some cases be expected to be possible to
improve these computational times in practice by tailoring parts of
the sparse linear algebra, but it would require additional coding.

\begin{figure}[tbh]
  \centering
  \includegraphics[width=0.44\textwidth]{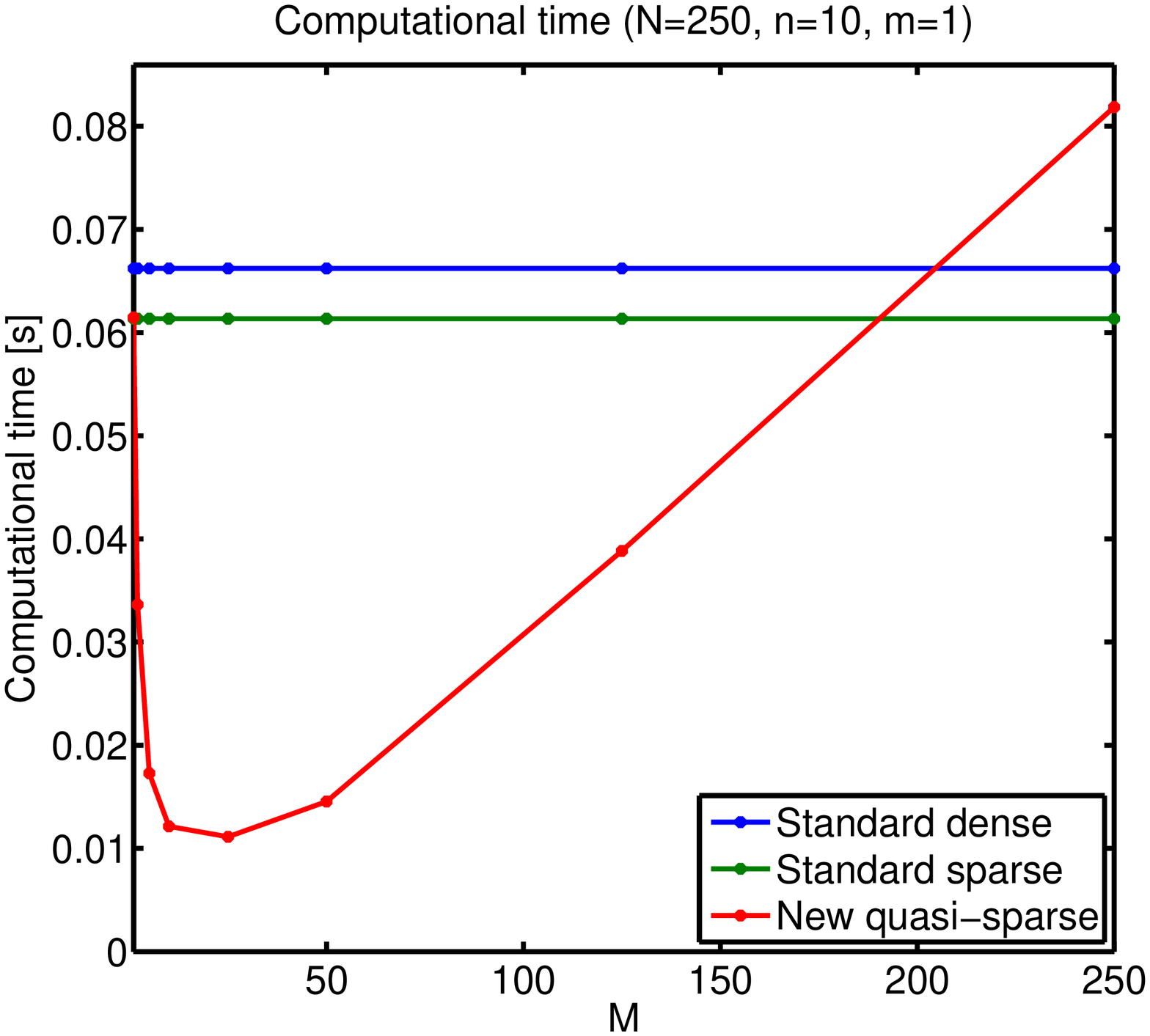}\hspace{0.05\textwidth}\includegraphics[width=0.45\textwidth]{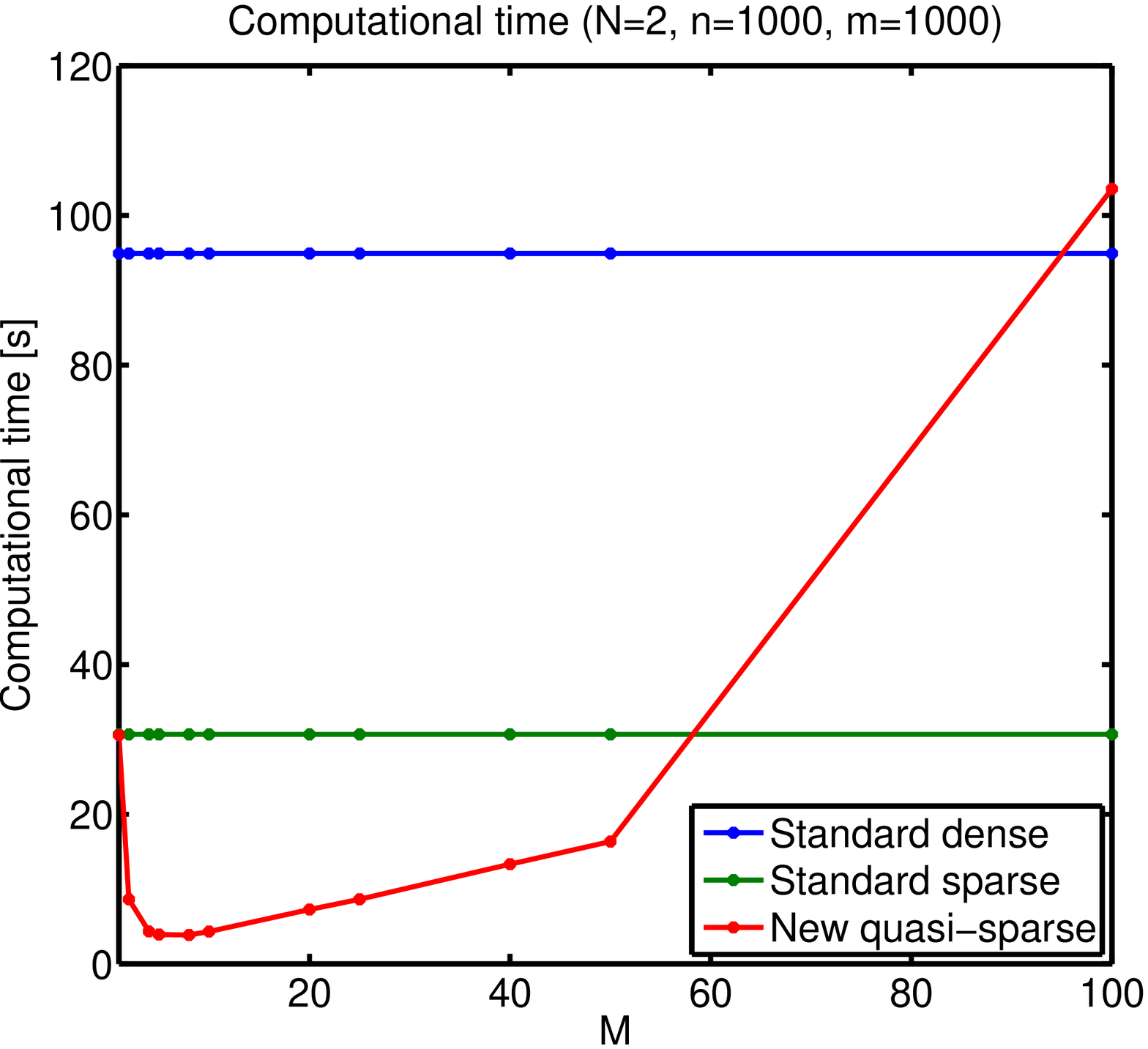}
  \caption{Computational time for the approach in
    Section~\ref{subsec:inc_blk_size} (left) and the approach in
    Section~\ref{subsec:dec_blk_size} (right) compared to standard
    dense and sparse approaches. The respective approaches presented
    in this work are in the plots denoted ``New quasi-sparse''.}
  \label{fig:quasi_sparse_comp_times}
\end{figure}

From the left plot in Figure~\ref{fig:quasi_sparse_comp_times}, it is
clear that the performance gain of the approach in
Section~\ref{subsec:inc_blk_size} can be significant for the chosen
example where $N=250$, $n=10$, and $m=1$. The chosen example is one in
which the traditional rule-of-thumb would have suggested that the
sparse approach would have been preferable. This is indeed the case,
but only with a small margin, and it turns out that the sparse and the
dense approaches require $61\,$ms and $66\,$ms, respectively. However,
by reformulating the problem using the results presented in
Section~\ref{subsec:inc_blk_size} and solving this new formulation
using the same sparse Riccati solver, the performance of the sparse
approach applied to the reformulated problem is significantly improved
and the computational time reduced to $11\,$ms for $M=25$.

From the right plot in Figure~\ref{fig:quasi_sparse_comp_times}, it is
clear that the performance gain of the approach in
Section~\ref{subsec:dec_blk_size} also can be significant for the
chosen example where $N=2$, $n=1000$, and $m=1000$. The computational
time for the standard sparse method applied to the reformulated
problem as proposed in Section~\ref{subsec:dec_blk_size} is $3.8\,$s
attained for $M=8$, which should be compared with the standard methods
that need $94\,$s and $31\,$s, respectively. Note also, that the
chosen example is one in which the traditional rule-of-thumb would
have suggested that the dense approach should have been
preferable. This is not the case and instead the performance of the
sparse one is actually better than the dense one.

\section{Conclusions}
In this work, two reformulations of MPC problems are presented. The
main idea in both reformulations is to trade-off the length of the
prediction horizon and the number of control signals in each step
along the horizon. This in turn affects the block-size used in the
numerical linear algebra and can as such a tool be used to find a
formulation of the problem that better utilizes the software libraries
and the hardware available on-line. It is shown that these new
formulations of the problem can significantly reduce the theoretical
number of flops and it is verified in numerical experiments that
significant reductions in computational times can be obtained also in
practice. More general, the result presented in this work shows that
the question to be answered when formulating an MPC problem is not
whether a sparse or a dense formulation should be used, but rather how
sparse the formulation should be.

\appendix

\section{Definition of stacked matrices}
\begin{equation}
  \label{eq:mat_defs}
  \begin{split}
    \timestack{x} &= 
    \begin{bmatrix} x_0^T & x_1^T & \hdots & x_N^T \end{bmatrix}^T, \;
    \timestack{u} = \begin{bmatrix} u_0^T & u_1^T & \hdots & u_{N-1}^T \end{bmatrix}^T \\
    \timestack{Q} &= \diag(Q_0,Q_1,\hdots,Q_N), \; \timestack{R} =
    \diag(R_0,R_1,\hdots,R_{N-1}) \\
    S_{x,k} &=
    \begin{bmatrix}
      I \\ A \\ A^2 \\ \vdots \\ A^k
    \end{bmatrix}, \;
    S_{u,k} =
    \begin{bmatrix}
      0 & 0 & \hdots & 0 \\
      B & 0 & \hdots & 0 \\
      AB & B & \hdots & 0 \\
      \vdots & \vdots & \ddots & \vdots \\
      A^{k-1}B & A^{k-2}B & \hdots & B
    \end{bmatrix}\\
    \timestack{A}&=
    \begin{bmatrix}
      -I & 0 & 0 & \hdots & 0 \\
      A & -I & 0 & \hdots & 0 \\
      0 & A & -I & \hdots & 0 \\
      0 & \ddots & \ddots & \ddots & 0 \\
      0 & \hdots & 0 & A & -I
    \end{bmatrix}
    , \;
    \timestack{A}^{-1}=
    \begin{bmatrix}
      -I & 0 & 0 & \hdots & 0 \\
      -A & -I & 0 & \hdots & 0 \\
      -A^2 & -A & -I & \hdots & 0 \\
      \vdots & \ddots & \ddots & \ddots & 0 \\
      -A^{N} & \hdots & 0 & -A & -I
    \end{bmatrix}
    , \\
    \timestack{B}&=
    \begin{bmatrix}
      0 & \hdots & \hdots &  0 \\
      B & 0 & \hdots & 0 \\
      0 & B & \hdots & 0 \\
      0 & \ddots & \ddots & 0 \\
      0 & \hdots & 0 & B
    \end{bmatrix}
    ,\\
    \timestack{b} &= \begin{bmatrix} -x^{0\,T} & 0 & 0 &
      \hdots & 0 \end{bmatrix}^T \\
    \timestack{H_x}&=\diag\parens{H_x(0),\hdots,H_x(N))}, \;
    \timestack{H_u}=\diag\parens{H_u(0),\hdots,H_u(N-1))}, \\
    \timestack{h}&=\brackets{h^T(0),\hdots,h^T(N)}^T
  \end{split}
\end{equation}
Note that, from~\eqref{eq:mat_defs} it follows that $S_{u,N} =
-\timestack{A}^{-1}\timestack{B}$.

\bibliography{IEEEabrv,references}

\end{document}